\documentstyle[11pt,twoside]{article}
\newfont{\Bbb}{msbm10 scaled\magstephalf}
\begin{document}
\title{The Essential Norm of Composition Operator between
Generalized Bloch  Spaces in Polydiscs and its Applications }
\author{Zehua Zhou and Yan Liu\\
Department of Mathematics, Tianjin University, Tianjin 300072, China\\
E-mail: zehuazhou2003@yahoo.com.cn and zhzhou@tju.edu.cn}
\date{}
\maketitle \footnote{Supported in part by National Natural Science
Foundation of China (Grand Nos. 10371091, 10001030), and LiuHui
Center for Applied Mathematics, Nankai University \& Tianjin
University}

\begin{abstract}Let $U^{n}$ be the unit polydisc of $\mbox{\Bbb C}^{n}$ and
$\phi=(\phi_1, \ldots, \phi_n)$ a holomorphic self-map of $U^{n}.$
By ${\cal B}^p(U^{n})$, ${\cal B}^p_{0}(U^{n})$ and ${\cal
B}^p_{0*}(U^{n})$ denote the $p$-Bloch space, Little $p$-Bloch
space and Little star $p$-Bloch space in the unit polydisc $U^n$
respectively, where $p, q>0$. This paper gives the estimates of
the essential norms of bounded composition operators $C_{\phi}$
induced by $\phi$ between ${\cal B}^p(U^n)$ (${\cal B}^p_{0}(U^n)$
or ${\cal B}^p_{0*}(U^n)$ ) and ${\cal B}^q(U^n)$ (${\cal
B}^q_{0}(U^n)$ or ${\cal B}^q_{0*}(U^n)$). As their applications,
some necessary and sufficient conditions for the bounded
composition operators $C_{\phi}$ to be compact from ${\cal
B}^p(U^n)$ $({\cal B}^p_{0}(U^n)$ or ${\cal B}^p_{0*}(U^n))$ into
${\cal B}^q(U^n)$ (${\cal B}^q_{0}(U^n)$ or ${\cal
B}^q_{0*}(U^n)$) are obtained.
\end{abstract}

\hspace{4mm}Keywords \hspace{2mm}Bloch space;
 Polydisc; Composition operator; Essential norm

\hspace{4mm}2000\hspace{1mm}Mathematics Subject Classification
\hspace{2mm} 47B38, 32A37, 47B33, 32A30

\section{Introduction}
\newtheorem{Definition}{Definition}
\newtheorem{Lemma}{Lemma}
\newtheorem{Theorem}{Theorem}
\newtheorem{Proposition}{Proposition}
\newtheorem{Corollary}{Corollary}

Let $\Omega$ be a bounded homogeneous domain in $\mbox{\Bbb
C}^{n}.$ The class of all holomorphic functions with domain
$\Omega$ will be denoted by $H(\Omega).$ Let $\phi $ be a
holomorphic self-map of $\Omega,$ the composition operator
$C_{\phi}$ induced by $\phi$ is defined by
$$(C_{\phi}f)(z)=f(\phi(z)),$$ for $z$ in $\Omega$ and $f\in
H(\Omega)$.

Let $K(z,z)$ be the Bergman kernel function of $\Omega$, the
Bergman metric $H_{z}(u,u)$ in $\Omega$ is defined by
$$H_{z}(u,u)=\displaystyle\frac{1}{2}
\sum\limits^{n}_{j,k=1}
\displaystyle\frac{\partial^{2}\log K(z,z)}{\partial z_{j}
\partial {\overline {z}}_{k}}u_{j}{\overline u}_{k},$$
where $z\in\Omega$
and $u=(u_{1},\ldots,u_{n})\in \mbox{\Bbb C}^{n}.$

Following Timoney [1], we say that $f\in H(\Omega)$ is in the
Bloch space ${\cal B}(\Omega),$ if
$$\|f\|_{{\cal B}(\Omega)}=\sup\limits_{z\in \Omega}Q_{f}(z)<\infty,$$
where
\begin{equation}Q_{f}(z)=\sup\left\{\displaystyle\frac
{|\bigtriangledown f(z)u|}{H^{\frac{1}{2}}_{z}(u,u)}: u\in
\mbox{\Bbb C}^{n}-\{0\}\right\},\label{1}\end{equation} and
$\bigtriangledown f(z) =\left(\frac{\partial f(z)}{\partial
z_{1}}, \ldots, \frac{\partial f(z)}{\partial z_{n}} \right),
\bigtriangledown f(z)u =\sum\limits^{n}_{l=1}\frac{\partial
f(z)}{\partial z_{l}}u_{l}.$

The little Bloch space ${\cal B}_0(\Omega)$ is the closure in the
Banach space ${\cal B}(\Omega)$ of the polynomial functions.

Let $\partial\Omega$ denote the boundary of $\Omega$. Following
Timoney [2], for $\Omega=B_n$ the unit ball of $\mbox{\Bbb C}^n$,
${\cal B}_0(B_n)=\left\{f\in{\cal B}(B_n): Q_f(z)\to 0,
\mbox{as}\hspace*{2mm}z\to\partial B_n \right\};$ for $\Omega=\cal
D$ the bounded symmetric domain other than the ball $B_n$,
$\left\{ f\in{\cal B}({\cal D}): Q_f(z)\to 0,
\mbox{as}\hspace*{2mm} z\to\partial{\cal D}\right\}$ is the set of
constant functions on $\cal D.$ So if $\cal D$ is a bounded
symmetric domain other than the ball, we denote the ${\cal
B}_{0*}({\cal D})= \left\{f\in{\cal B}({\cal D}): Q_f(z)\to 0,
\mbox{as}\hspace*{2mm} z\to\partial^*{\cal D}\right\}$ and call it
little star Bloch space, here $\partial^*{\cal D}$ means the
distinguished boundary of $\cal D$. The unit ball is the only
bounded symmetric domain $\cal D$ with the property that
$\partial^*{\cal D}=\partial{\cal D}.$

Let $U^n$ be the unit polydisc of $\mbox{\Bbb C}^n$. Timoney [1]
shows that $f\in{\cal B}(U^n)$ if and only if
$$\|f\|_1=|f(0)|+\sup\limits_{z\in U^n}\sum\limits^n_{k=1}
\left|\frac{\partial f} {\partial z_k}(z)\right|(1-
|z_k|^2)<+\infty,$$ where $f\in H(U^n).$

This definition was the starting point for introducing the
$p$-Bloch spaces.

Let $p>0,$ a function $f\in H(U^n)$ is said to belong to the
$p$-Bloch space ${\cal B}^p(U^n)$ if
$$\|f\|_p=|f(0)|+\sup\limits_{z\in U^n}\sum\limits^n_{k=1}
\left|\frac{\partial f} {\partial z_k}(z)\right|\left(1-
|z_k|^2\right)^p<+\infty.$$ It is easy to show that ${\cal
B}^p(U^n)$ is a Banach space with the norm $\|\cdot\|_p.$

Just like Timoney [2], if
$$\lim_{z\to\partial U^n}\sum\limits^n_{k=1}\left|\frac{\partial f}{\partial
z_k}(z)\right|(1- |z_k|^2)^p=0,$$ it is easy to show that $f$ must
be a constant. Indeed, for fixed $z_1\in U,$
$\displaystyle\frac{\partial f}{\partial z_1}(z)(1-|z_1|^2)^p$ is
a holomorphic function in $z'=(z_2,\cdots,z_n)\in U^{n-1}$. If
$z\to\partial U^n$, then $z'\to\partial U^{n-1},$ which implies
that
$$\lim\limits_{z'\to\partial U^{n-1}}\left|\frac{\partial f} {\partial
z_1}(z)\right|\left(1- |z_1|^2\right)^p=0.$$ Hence,
$\frac{\partial f} {\partial z_1}(z)\left(1-
|z_1|^2\right)^p\equiv 0$ for every $z'\in \partial U^{n-1},$ and
for each $z_1\in U,$ and consequently $\frac{\partial f} {\partial
z_1}(z)=0$ for every $z\in U^n.$ Similarly, we can obtain that
$\frac{\partial f} {\partial z_j}(z)=0$ for every $z_j\in U^n$ and
each $j\in\{2,\cdots,n\},$ therefore $f\equiv const .$

So, there is no sense to introduce the corresponding little
$p$-Bolch space in this way. We will say that the little $p$-Bolch
space ${\cal B}_0^p(U^n)$ is the closure of the polynomials in the
$p$-Bolch space. If $f\in H(U^n)$ and
$$\sup\limits_{z\in \partial^*U^n}\sum\limits^n_{k=1}
\left|\displaystyle\frac{\partial f} {\partial
z_k}(z)\right|\left(1- |z_k|^2\right)^p=0,$$ we say $f$ belongs to
little star $p$-Bolch space ${\cal B}_{0*}^p(U^n).$ Using the same
methods as that of Theorem  4.14 in reference [2], we can show
that ${\cal B}^p_{0}(U^n)$ is a proper subspace of ${\cal
B}^p_{0*}(U^n)$ and ${\cal B}^p_{0*}(U^n)$ is a non-separable
closed subspace of ${\cal B}^p(U^n).$

Let $\phi $ be a holomorphic self-map of $U^n,$ the composition
operator $C_{\phi}$ induced by $\phi$ is defined by
$(C_{\phi}f)(z)=f(\phi(z))$ for $z$ in $U^n$ and $f\in H(U^n)$.
For the unit disc $U\subset\mbox{\Bbb C},$ Madigan and Matheson
[3] proved that $C_{\phi}$ is always bounded on ${\cal B}(U)$ and
bounded on ${\cal B}_0(U)$ if and only if $\phi\in{\cal B}_0(U).$
They also gave the sufficient and necessary conditions that
$C_{\phi}$ is compact on ${\cal B}(U)$ or ${\cal B}_0(U).$ More
recently, [4,5,7] gave some sufficient and necessary conditions
for $C_{\phi}$ to be compact on the Bloch spaces in polydisc.

We recall that the essential norm of a continuous linear operator $T$ is
the distance from $T$ to the compact operators, that is,
\begin{equation}\|T\|_e=\inf\{\|T-K\|: K \mbox{ is compact}\}.\label{2}\end{equation}
Notice that $\|T\|_e=0$ if and only if $T$ is compact, so that estimates
on $\|T\|_e$ lead to conditions for $T$ to be compact.

In this paper, we give some estimates of the essential norms of
bounded composition operators $C_{\phi}$ between ${\cal B}^p(U^n)$
(${\cal B}^p_{0}(U^n)$ or ${\cal B}^p_{0*}(U^n)$ ) and ${\cal
B}^q(U^n)$ (${\cal B}^q_{0}(U^n)$ or ${\cal B}^q_{0*}(U^n)$). As
their consequences, some necessary and sufficient conditions for
the bounded composition operators $C_{\phi}$ to be compact from
${\cal B}^p(U^n)$ (${\cal B}^p_{0}(U^n)$ or ${\cal B}^p_{0*}(U^n)$
) into ${\cal B}^q(U^n)$ (${\cal B}^q_{0}(U^n)$ or ${\cal
B}^q_{0*}(U^n)$) are obtained.

The fundamental ideals of the proof are those used by J. H.
Shpairo [8] to obtain the essential norm of a composition operator
on Hilbert spaces of analytic functions (Hardy and weighted
Bergman spaces) in terms of natural counting functions associated
with $\phi$. This paper generalizes the result on the Bloch space
in [10] to the Bloch-type space in polydisk.

Throughout the remainder of this paper $C$ will denote a positive
constant, the exact value of which will vary from one appearance
to the next.

Our main results are the following:

\begin{Theorem} Let $\phi=(\phi_1, \phi_2, \cdots, \phi_n)$ be a holomorphic
self-map of $U^n$ and $\|C_{\phi}\|_e$ the essential norm of a
bounded composition operator $C_{\phi}:$ ${\cal B}^p(U^n)$ (
${\cal B}^p_{0}(U^n)$ or ${\cal B}^p_{0*}(U^n)$)$\rightarrow$ $
{\cal B}^q(U^n)$ ( ${\cal B}^q_{0}(U^n)$ or ${\cal
B}^q_{0*}(U^n)$) , then
\begin{eqnarray}&&\displaystyle\frac{1}{n}\lim\limits_{\delta\to 0}
\sup\limits_{dist(\phi(z),\partial U^n)<\delta}\sum\limits^n_{k,l=1}
\left|\displaystyle\frac{\partial \phi_{l}}{\partial z_k}(z)\right|
\displaystyle\frac{(1-|z_k|^2)^q}{(1-|\phi_l(z)|^2)^p}\nonumber\\
&&\leq\|C_{\phi}\|_e \leq 2\lim\limits_{\delta\to 0}
\sup\limits_{dist(\phi(z),\partial
U^n)<\delta}\sum\limits^n_{k,l=1}
\left|\displaystyle\frac{\partial \phi_{l}}{\partial
z_k}(z)\right|
\displaystyle\frac{(1-|z_k|^2)^q}{(1-|\phi_l(z)|^2)^p}.\label{3}\end{eqnarray}
\end{Theorem}

By Theorem 1 and the fact that $C_{\phi}:$ ${\cal B}^p(U^n)$ (or
${\cal B}^p_{0}(U^n)$ or ${\cal B}^p_{0*}(U^n)$)$\rightarrow$ $
{\cal B}^q(U^n)$ (or ${\cal B}^q_{0}(U^n)$ or ${\cal
B}^q_{0*}(U^n)$) is compact if and only if $\|C_{\phi}\|_e=0$, we
obtain Theorem 2 at once.

\begin{Theorem}\hspace{2mm}Let $\phi=(\phi_1, \ldots, \phi_n)$ be a
holomorphic self-map of $U^{n}.$ Then the bounded composition
operator $C_{\phi}:$ ${\cal B}^p(U^n)$ (${\cal B}^p_{0}(U^n)$ or
${\cal B}^p_{0*}(U^n)$)$\rightarrow$ $ {\cal B}^q(U^n)$ (${\cal
B}^q_{0}(U^n)$ or ${\cal B}^q_{0*}(U^n)$) is compact if and only
if for any $\varepsilon>0,$ there exists a $\delta$ with
$0<\delta<1,$ such that
\begin{equation}\sup\limits_{dist(\phi(z),\partial U^n)<\delta}
\sum\limits^n_{k,l=1} \left|\displaystyle\frac{\partial
\phi_{l}}{\partial z_k}(z)\right|
\displaystyle\frac{(1-|z_k|^2)^q}{(1-|\phi_l(z)|^2)^p}<\varepsilon.\label{4}\end{equation}
\end{Theorem}

When $n=1,$ on ${\cal B}(U)$ we obtain Theorem 2 in [3]. Since
$\partial U=\partial^* U,$ ${\cal B}_0(U)={\cal B}_{0*}(U),$ we
can also obtain Theorem 1 in [3].

By Theorem 2 and Lemmas 3, 4 and 5 in next part, we can get the
following three Corollaries.

\begin{Corollary} Let $\phi=(\phi_1, \ldots, \phi_n)$ be a
holomorphic self-map of $U^{n}.$ Then\\ $C_{\phi}:{\cal
B}^p(U^{n})$(${\cal B}^p_{0}(U^n)$ or ${\cal
B}^p_{0*}(U^n)$)$\rightarrow{\cal B}^q(U^{n})$ is compact if and
only if
$$\sum\limits^n_{k,l=1} \left|\displaystyle\frac{\partial
\phi_{l}}{\partial z_k}(z)\right|
\displaystyle\frac{(1-|z_k|^2)^q}{(1-|\phi_l(z)|^2)^p}\leq C$$ for
all $z\in U^n$ and (\ref{4}) holds.\end{Corollary}

{\bf Proof}\hspace*{4mm} By Lemma 3 in next part, we know
$C_{\phi}:{\cal B}^p(U^{n})$(${\cal B}^p_{0}(U^n)$ or ${\cal
B}^p_{0*}(U^n)$)$\rightarrow{\cal B}^q(U^{n})$ is bounded. It
follows from Theorem 2 that $C_{\phi}:{\cal B}^p(U^{n})$(${\cal
B}^p_{0}(U^n)$ or ${\cal B}^p_{0*}(U^n)$)$\rightarrow{\cal
B}^q(U^{n})$ is compact.

Conversely, if $C_{\phi}:{\cal B}^p(U^{n})$(${\cal B}^p_{0}(U^n)$
or ${\cal B}^p_{0*}(U^n)$)$\rightarrow{\cal B}^q(U^{n})$ is
compact, it is clear that $C_{\phi}:{\cal B}^p(U^{n})$(${\cal
B}^p_{0}(U^n)$ or ${\cal B}^p_{0*}(U^n)$)$\rightarrow{\cal
B}^q(U^{n})$ is bounded, by Theorem 2, (\ref{4}) holds.

\begin{Corollary} Let $\phi=(\phi_1, \ldots, \phi_n)$ be a
holomorphic self-map of $U^{n}.$ Then \\
$C_{\phi}:$${\cal B}^p_{0*}(U^{n})$(${\cal
B}^p_{0}(U^n)$)$\rightarrow{\cal B}^q_{0*}(U^{n})$ is compact if
and only if $\phi_l\in {\cal B}^q_{0*}(U^n)$ for every
$l=1,2,\cdots, n$ and (\ref{4}) holds.\end{Corollary}

{\bf Proof}\hspace*{4mm} Note that Lemma 4 in next part, similar
to the proof of Corollary 1, the Corollary follows.

\begin{Corollary} Let $\phi=(\phi_1, \ldots, \phi_n)$ be a
holomorphic self-map of $U^{n}.$ Then \\
$C_{\phi}:$${\cal B}^p_{0}(U^{n})\rightarrow{\cal B}^q_{0}(U^{n})$
is compact if and only if $\phi_l\in {\cal B}^q_{0}(U^n)$ for
every $l=1,2,\cdots, n$ and (\ref{4}) holds.\end{Corollary}

{\bf Proof}\hspace*{4mm} Note that Lemma 5 in next part, similar
to the proof of Corollary 1, the Corollary follows.

\section{Some Lemmas}
In order to prove Theorem 1, we need some Lemmas.
\begin{Lemma} Let $f\in{\cal B}^p(U^n),$ then

(1) If $0\leq p<1,$ then $\|f(z)|\leq
|f(0)|+\displaystyle\frac{n}{1-p}\|f\|_p;$

(2) If $p=1,$ then $|f(z)|\leq \left(1+\displaystyle\frac{1}{n\ln
2}\right)\left(\sum\limits^n_{k=1}\ln
\displaystyle\frac{2}{1-|z_k|^2}\right)\|f\|_p.$

(3) If $p>1,$ then $|f(z)|\leq
\left(\displaystyle\frac{1}{n}+\displaystyle\frac{2^{p-1}}{p-1}\right)
\sum\limits^n_{k=1}\displaystyle\frac{1}{(1-|z_k|^2)^{p-1}}\|f\|_p.$
\end{Lemma}

{\bf Proof}\hspace{2mm} This Lemma can be proved by some integral
estimates (if necessary, the proof can be omitted).

By the definition of $\|.\|_{p}$,
$$|f(0)|\leq \|f\|_{p},\hspace*{4mm}\left|\displaystyle\frac{\partial f(z)
}{\partial z_l}\right|\leq
\displaystyle\frac{\|f\|_{p}}{(1-|z_l|^2)^p}
\hspace*{4mm}(l\in\{1,2,\cdots,n\})$$ and
\begin{eqnarray*}&&f(z)-f(0)=
\int^1_0 \displaystyle\frac{d f(tz)}{d
t}dt=\sum\limits^n_{l=1}\int^1_0 z_l\displaystyle\frac{\partial
f}{\partial\zeta_l}(tz)dt,\end{eqnarray*}
So\begin{eqnarray}&&|f(z)|\leq
|f(0)|+\sum\limits^n_{l=1}|z_l|\int^1_0\displaystyle\frac{\|f\|_p}{\left(1-t^2|z_l|^2\right)^p}dt
\nonumber\\
&&\leq
\|f\|_p+\|f\|_p\sum\limits^n_{l=1}\int^{|z_l|}_0\displaystyle\frac{1}{\left(1-t^2\right)^p}dt.\label{5}
\end{eqnarray}
If $p=1,$
\begin{equation}\int^{|z_l|}_0\displaystyle\frac{1}{\left(1-t^2\right)^p}dt=\displaystyle\frac{1}{2}
\ln\displaystyle\frac{1+|z_l|}{1-|z_l|}\leq
\displaystyle\frac{1}{2}
\ln\displaystyle\frac{4}{1-|z_l|^2}.\label{6}\end{equation} It is
clear that $\ln\displaystyle\frac{4}{1-|z_l|^2}>\ln 4=2\ln 2,$
so\begin{equation}1\leq \displaystyle\frac{1}{2\ln 2}
\ln\displaystyle\frac{4}{1-|z_l|^2}\leq
\displaystyle\frac{1}{2n\ln 2}
\sum\limits^n_{l=1}\ln\displaystyle\frac{4}{1-|z_l|^2}.\label{7}\end{equation}
Combining (\ref{5}),(\ref{6}) and (\ref{7}), we get
$$|f(z)|\leq \left(\displaystyle\frac{1}{2}+\displaystyle\frac{1}{2n\ln
2}\right)\left(\sum\limits^n_{l=1}\ln\displaystyle\frac{4}{1-|z_l|^2}\right)\|f\|_p.$$
If $p\neq 1,$
\begin{eqnarray}&&\int^{|z_l|}_0\displaystyle\frac{1}{\left(1-t^2\right)^p}dt
=\int^{|z_l|}_0\displaystyle\frac{1}{(1-t)^p}\cdot
\displaystyle\frac{1}{(1+t)^p}dt\nonumber\\
&&\leq
\int^{|z_l|}_0\displaystyle\frac{1}{(1-t)^p}dt=\displaystyle\frac{1-(1-|z_l|)^{-p+1}}{1-p}.\label{8}
\end{eqnarray}
If $0<p<1,$ (\ref{8}) gives that
$\int^{|z_l|}_0\displaystyle\frac{1}{\left(1-t^2\right)^p}dt\leq\displaystyle\frac{1}{1-p},$
it follows from (\ref{5}) that
$|f(z)|\leq\left(1+\displaystyle\frac{n}{1-p}\right)\|f\|_p.$

If $p>1,$ (\ref{8}) gives that
$$\int^{|z_l|}_0\displaystyle\frac{1}{\left(1-t^2\right)^p}dt\leq \displaystyle\frac{1-(1-|z_l|^{p-1})}{(p-1)(1-|z_l|)^{p-1}}
\leq\displaystyle\frac{2^{p-1}}{(p-1)(1-|z_l^2|)^{p-1}},$$ it
follows from (\ref{5}) that
\begin{eqnarray*}|f(z)|&\leq&\|f\|_p+\displaystyle\frac{2^{p-1}}{p-1}\left(\sum\limits^n_{l=1}
\displaystyle\frac{1}{(1-|z_l|^2)^{p-1}}\right)\|f\|_p\\
&\leq&
\left(\displaystyle\frac{1}{n}+\displaystyle\frac{2^{p-1}}{p-1}\right)\left(\sum\limits^n_{l=1}
\displaystyle\frac{1}{(1-|z_l|^2)^{p-1}}\right)
\|f\|_{p}.\end{eqnarray*}Now the Lemma is proved.

\begin{Lemma} Set $$f_w(z)=\int_0^{z_l}\frac {dt}{(1-\bar wt)^p},$$ where $w\in
U.$  Then $f\in {\cal B}^p_0(U^n)\subset{\cal
B}^p_{0*}(U^n)\subset {\cal B}^p(U^n).$ \end{Lemma}

{\bf Proof}\hspace{2mm}Since
$$\displaystyle\frac{\partial
f_w}{\partial z_l} =\left(1-\overline wt\right)^{-p},
\hspace*{4mm} \displaystyle\frac{\partial f_w}{\partial z_i}=0,
\hspace*{4mm}(i\neq l),$$ it follows that
$$|f(0)|+\sum\limits^n_{k=1}\left|\displaystyle\frac{\partial f_w}{\partial z_k}(z)\right|(1-|z_k|^2)^p
=\displaystyle\frac{(1-|z_l|^2)^p}{|1-\overline
wt|^p}\leq(1+|z_l|^2)^p\leq 2^p.$$ Hence $f_w\in {\cal B}^p(U^n).$

Now we prove that $f_w\in{\cal B}^p_0(U^n).$ Using the asymptotic
formula
$$(1-\bar w t)^{-p}=\sum\limits^{+\infty}_{k=0}\frac{p(p+1)\cdots
(p+k-1)}{k!}(\bar w)^kt^k,$$ we obtain
$$f_w(z)=\sum\limits^{+\infty}_{k=0}\frac{p(p+1)\cdots
(p+k-1)}{k!}(\bar w)^k\int^{z_l}_0 t^kdt.$$ Denote
$P_n(z)=\sum\limits^{n}_{k=0}\frac{p(p+1)\cdots (p+k-1)}{k!}(\bar
w)^k\int^{z_l}_0t^kdt,$  it is easy to see that
$$f_w(z)-P_n(z)=\sum\limits^{+\infty}_{k=n+1}\frac{p(p+1)\cdots (p+k-1)}{k!}(\bar w)^k\int^{z_l}_0t^k dt,$$
$$\left|\frac{\partial (f_w-P_n)}{\partial z_l}\right|\leq
\sum\limits^{+\infty}_{k=n+1}\frac{p(p+1)\cdots (p+k-1)}{k!}|w|^k
\to 0, \mbox{as}\ \  n\to\infty,$$
\begin{eqnarray*}\|f_w-P_n\|_p&=&|f_w(0)-P_n(0)|+\sup\limits_{z\in
U^n}\left|\frac{\partial (f_w-P_n)}{\partial
z_l}\right|(1-|z_l|^2)^p \\
&\leq& \sup\limits_{z\in U^n}\left|\frac{\partial
(f_w-P_n)}{\partial z_l}\right|\to 0, \end{eqnarray*} it shows
that $f_w\in {\cal B}^p_0(U^n).$ So $f\in {\cal
B}^p_0(U^n)\subset{\cal B}^p_{0*}(U^n)\subset {\cal B}^p(U^n).$

\begin{Lemma} Let $\phi=(\phi_1, \ldots, \phi_n)$ be a holomorphic self-map
of $U^n$, $p,q>0.$ Then $C_{\phi}: {\cal B}^p(U^n) ({\cal
B}^p_0(U^n)$ or ${\cal B}^p_{0*}(U^n))\longrightarrow {\cal
B}^q(U^n) $ is bounded if and only if there exists a constant $C$
such that
\begin{equation}\sum\limits^n_{k,l=1}\left|\displaystyle\frac{\partial \phi_{l}}
{\partial
z_k}(z)\right|\displaystyle\frac{\left(1-|z_k|^2\right)^q}
{\left(1-|\phi_l(z)|^2\right)^p}\leq C ,\label{9}\end{equation}for
all $z\in U^n.$\end{Lemma}

{\bf Proof}\hspace*{4mm}First assume that condition (\ref{9})
holds. Let $f\in{\cal B}^p(U^n)({\cal B}^p_0(U^n)$ or ${\cal
B}^p_{0*}(U^n)),$ by Lemma 1, we know the evaluation at $\phi(0)$
is a bounded linear functional on ${\cal B}^p(U^n),$ so
$|f(\phi(0))|\leq C\|f\|_p.$

On the other hand we have
\begin{eqnarray}&&\sum\limits^n_{k=1}\left|\frac{\partial
\left(C_{\phi}f(z)\right)} {\partial z_k}\right| (1-|z_k|^2)^q
=\sum\limits^n_{k=1}\left|\sum\limits^n_{l=1}\frac{\partial
f}{\partial \phi_l}(\phi(z))
\frac{\partial \phi_l}{\partial z_k}(z)\right|(1-|z_k|^2)^q\nonumber\\
&&\leq \sum\limits^n_{k, l=1}\left|\frac{\partial f}{\partial
\phi_l}(\phi(z))\frac{\partial
\phi_l}{\partial z_k}(z)\right|(1-|z_k|^2)^q\nonumber\\
&&\leq \sum\limits^n_{l=1}\left|\frac{\partial f}{\partial
\phi_l}(\phi(z))\right|\left(1-|\phi_l(z)|^2\right)^p
\sum\limits^n_{k,l=1}\left|\frac{\partial \phi_{l}} {\partial
z_k}(z)\right|\frac{\left(1-|z_k|^2\right)^q}
{\left(1-|\phi_l(z)|^2\right)^p}\label{10}\\
&&\leq\|f\|_p\sum\limits^n_{k,l=1}\left|\frac{\partial \phi_l}
{\partial z_k}(z)\right|
\frac{\left(1-|z_k|^2\right)^q}{\left(1-|\phi_l(z)|^2\right)^p}.\label{11}
\end{eqnarray}

From (\ref{11}) it follows that
$$\sum\limits^n_{k=1}\left|\frac{\partial
\left(C_{\phi}f(z)\right)} {\partial z_k}\right| (1-|z_k|^2)^q\leq
C\|f\|_p.$$ So $C_{\phi}: {\cal B}^p(U^n)\to {\cal B}^q(U^n) $ is
bounded.

For the converse, assume that $C_{\phi}: {\cal B}^p(U^n) ({\cal
B}^p_0(U^n)$ or ${\cal B}^p_{0*}(U^n))\longrightarrow {\cal
B}^q(U^n) $ is bounded, with
\begin{equation}\|C_{\phi}f\|_q\leq C\|f\|_p\label{12}\end{equation}
for all $f\in{\cal B}^p(U^n)({\cal B}^p_0(U^n)$ or ${\cal
B}^p_{0*}(U^n)).$

For fixed $l (1\leq l\leq n),$ we will make use of a family of
test functions $\{f_{w}: w\in\mbox{\Bbb C}, |w|<1\}$ in ${\cal
B}(U^n)$ defined as follows: If $p>0$, let
$$f_w(z)=\int^{z_l}_0\left(1-\overline wz_l\right)^{-p}dt.$$
It follows from Lemma 2 that $$f_w\in {\cal
B}^p_0(U^n)\subset{\cal B}^p_{0*}(U^n)\subset {\cal B}^p(U^n).$$

For $z\in U^n,$ it follows from (\ref{12}) that
\begin{equation}\sum\limits^n_{k=1}\left|\sum\limits^n_{l=1}\displaystyle\frac{\partial
f_w(\phi(z))}{\partial \phi_l} \displaystyle\frac{\partial
\phi_l}{\partial z_k}(z)\right|(1-|z_k|^2)^q\leq
C.\label{13}\end{equation}

Let $w=\phi_l(z),$ then
$$\sum\limits^n_{k=1}\left|\displaystyle\frac{\partial
\phi_{l}} {\partial z_k}(z)\right|
\displaystyle\frac{\left(1-|z_k|^2\right)^q}{\left(1-|\phi_l(z)|^2\right)^p}
\leq C.$$ Now the proof of Lemma 3 is completed.

\begin{Lemma} Let $\phi=(\phi_1, \phi_2, \cdots, \phi_n)$ be a holomorphic
self-map of $U^n.$ Then \\
$C_{\phi}:{\cal B}^p_{0*}(U^{n}) ({\cal B}^p_0(U^n))
\rightarrow{\cal B}^q_{0*}(U^{n})$ is bounded if and only if
$\phi_l\in{\cal B}^q_{0*}(U^n)$ for every $l=1,2,\cdots, n$ and
(\ref{9}) holds.\end{Lemma}

{\bf Proof}\hspace*{4mm}If $C_{\phi}:{\cal B}^p_{0*}(U^{n}) ({\cal
B}^p_0(U^n)) \rightarrow {\cal B}^q_{0*}(U^{n})$  is bounded , it
is clear that, for every $l=1,2,\cdots, n$, $f_l(z)=z_l\in{\cal
B}^p_0(U^n)\subset{\cal B}^q_{0*}(U^n),$ so
$C_{\phi}f_l=\phi_l\in{\cal B}^q_{0*}(U^n).$ In the proof of Lemma
3, note that the test functions $f_w\in {\cal
B}^p_0(U^n)\subset{\cal B}^p_{0*}(U^n),$ we know $(\ref{9})$
holds.

In order to prove the Converse, we first prove that if
$\phi_l\in{\cal B}^q_{0*}(U^n)$ for every $l=1,2,\cdots,n.,$ then
$f\circ\phi\in{\cal B}^q_{0*}(U^n)$ for any $f\in{\cal
B}^p_{0*}(U^n).$

Without loss of generality, we prove this result when $n=2.$

For any sequence $\{z^j=(z^j_1, z^j_2)\}\subset U^n$ with
$z^j\to\partial^* U^n$ as $j\to\infty,$ then $$|z^j_1|\to 1,
|z^j_2|\to 1.$$ Since $|\phi_1(z^j)|<1$ and $|\phi_2(z^j)|<1,$
there exists a subsequence $\{z^{j_s}\}$ in $\{z^j\}$  such that
$$|\phi_1(z^{j_s})|\to \rho_1, |\phi_2(z^{j_s})|\to\rho_2,$$ as
$s\to\infty .$

It is clear that $0\leq\rho_1, \rho_2\leq 1.$
\begin{eqnarray}&&\left|\displaystyle\frac{\partial(f\circ\phi)}
{\partial z_k}(z^{j_s})\right|(1-|z^{j_s}_k|^2)^q\nonumber\\
&&\leq\left| \displaystyle\frac{\partial f}{\partial
w_1}(\phi(z^{j_s}))\right| \left|\displaystyle\frac
{\partial\phi_1}{\partial
z_k}(z^{j_s})\right|(1-|z^{j_s}_k|^2)^q+\left|
\displaystyle\frac{\partial f}{\partial w_2}(\phi(z^{j_s}))\right|
\left|\displaystyle\frac
{\partial\phi_2}{\partial z_k}(z^{j_s})\right|(1-|z^{j_s}_k|^2)^q\nonumber\\
&&=\left| \displaystyle\frac{\partial f}{\partial
w_1}(\phi(z^{j_s})) \right|(1-|\phi_1(z^{j_s})|^2)^p
\left|\displaystyle\frac {\partial\phi_1}{\partial
z_k}(z^{j_s})\right|\displaystyle\frac{(1-|z^{j_s}_k|^2)^q}
{(1-|\phi_1(z^{j_s})|^2)^p}\nonumber\\
&&+\left| \displaystyle\frac{\partial f}{\partial
w_2}(\phi(z^{j_s}))\right| (1-|\phi_2(z^{j_s})|^2)^p
\left|\displaystyle\frac {\partial\phi_2}{\partial
z_k}(z^{j_s})\right|\displaystyle\frac{(1-|z^{j_s}_k|^2)^q}
{(1-|\phi_2(z^{j_s})|^2)^p},\label{14}\end{eqnarray} $k=1,2.$

Now we prove the left of $(\ref{14})\to 0$ as $s\to\infty$
according to four cases.

Case 1. If $\rho_1<1$ and $\rho_2<1.$ It is clear that there exist
$r_1$ and $r_2$ such that $\rho_1<r_1<1$ and $\rho_2<r_2<1,$ so as
$j$ is large enough, $|\phi_1(z^{j_s})|\leq r_1$ and
$|\phi_2(z^{j_s})|\leq r_2.$

By $\phi_1, \phi_2\in{\cal B}^q_{0*}(U^n)$ and (\ref{14}), we get
\begin{eqnarray*}\left|\displaystyle\frac{\partial(f\circ\phi)}
{\partial z_k}(z^{j_s})\right|(1-|z^{j_s}_k|^2)^q &\leq& \|f\|_p
\displaystyle\frac{1}{(1-r_1^2)^p}\left|\displaystyle\frac
{\partial\phi_1}{\partial z_k}(z^{j_s})\right|(1-|z^{j_s}_k|^2)^q\\
&&+\|f\|_p\displaystyle\frac{1}{(1-r_2^2)^p}\left|\displaystyle\frac
{\partial\phi_2}{\partial z_k}(z^{j_s})\right|(1-|z^{j_s}_k|^2)^q\\
& &\to 0\end{eqnarray*} as $s\to\infty.$

Case 2. If $\rho_1=1$ and $\rho_2=1.$ Then
$\phi(z^{j_s})\to\partial^*U^n,$ by (\ref{9}) and $f\in{\cal
B}^p_{0*}(U^n)$, (\ref{14}) gives that
\begin{eqnarray*}&&C\left|\displaystyle\frac{\partial(f\circ\phi)} {\partial
z_k}(z^{j_s})\right|(1-|z^{j_s}_k|^2)^q\\
&&C\leq \left| \displaystyle\frac{\partial f}{\partial
w_1}(\phi(z^{j_s})) \right|(1-|\phi_1(z^{j_s})|^2)^p +\left|
\displaystyle\frac{\partial f}{\partial w_2}(\phi(z^{j_s}))\right|
(1-|\phi_2(z^{j_s})|^2)^p\to 0\end{eqnarray*} as $s\to\infty.$

Case 3. If $\rho_1<1$ and $\rho_2=1.$ Similar to Case 1,  we can
prove that
\begin{eqnarray}&&\left|
\displaystyle\frac{\partial f}{\partial w_1}(\phi(z^{j_s}))
\right|(1-|\phi_1(z^{j_s})|^2)^p\left|\displaystyle\frac
{\partial\phi_1}{\partial
z_k}(z^{j_s})\right|\displaystyle\frac{(1-|z^{j_s}_k|^2)^q}
{(1-|\phi_1(z^{j_s})|^2)^p}\nonumber\\
&&\leq\|f\|_p\displaystyle\frac{1}{(1-r_1^2)^p}
\left|\displaystyle\frac {\partial\phi_1}{\partial
z_k}(z^{j_s})\right|\displaystyle\frac{(1-|z^{j_s}_k|^2)^q}
{(1-|\phi_1(z^{j_s})|^2)^p} \to 0 \label{15}\end{eqnarray} as
$s\to\infty.$

On the other hand, for fixed $s,$ let $w^{j_s}_2=\phi_2(z^{j_s}),$
then $|w^{j_s}_2|<1.$ Denote
$$F(w_1)=\displaystyle\frac{\partial f}{\partial w_2}(w_1, w^{j_s}_2).$$
It is clear that $F(w_1)$ is holomorphic on $|w_1|<1,$ choose
$R_{j_s}\to 1$ with $r_1\leq R_{j_s}<1.$ $|\phi_1(z^{j_s})|\leq
r_1,$ so
$$|F(\phi_1(z^{j_s}))|\leq\max\limits_{|w_1|\leq r_1}|F(w_1)|
\leq\max\limits_{|w_1|\leq R_{j_s}}|F(w_1)|=
\max\limits_{|w_1|=R_{j_s}}|F(w_1)|=|F(w^{j_s}_1)|,$$ where
$|w^{j_s}_1|=R_{j_s}\to 1.$ This means that
$\left|\displaystyle\frac{\partial f}{\partial
w_2}(\phi_1(z^{j_s}), \phi_2(z^{j_s})) \right|\leq
\left|\displaystyle\frac{\partial f}{\partial w_2}(w^{j_s}_1,
w^{j_s}_2) \right|. $ Since $|w^{j_s}_1|\to 1,
|w^{j_s}_2|\to\rho_2=1$ and $f\in{\cal B}^p_{0*}(U^n),$
$$\left|\displaystyle\frac{\partial f}{\partial w_2}(w^{j_s}_1, w^{j_s}_2)\right|
(1-|w^{j_s}_2|^2)^p\to 0$$ as $s\to\infty,$ so by (\ref{9}),
\begin{eqnarray}&&\left|
\displaystyle\frac{\partial f}{\partial w_2}(\phi(z^{j_s}))\right|
(1-|\phi_2(z^{j_s})|^2) ^p\left|\displaystyle\frac
{\partial\phi_2}{\partial
z_k}(z^{j_s})\right|\displaystyle\frac{(1-|z^{j_s}_k|^2)^q}
{(1-|\phi_2(z^{j_s})|^2)^p}\nonumber\\
&&\leq C\left| \displaystyle\frac{\partial f}{\partial
w_2}(w^{j_s}_1, w^{j_s}_2)\right| (1-|w^{j_s}_2|^2)^p\to
0\label{16}\end{eqnarray} as $s\to\infty.$

By (\ref{15}) and (\ref{16}), (\ref{14}) gives
$$\left|\displaystyle\frac{\partial(f\circ\phi)}
{\partial z_k}(z^{j_s})\right|(1-|z^{j_s}_k|^2)^q\to 0,$$ as
$s\to\infty.$

Case 4. If $\rho_1=1$ and $\rho_2<1.$ Similar to Case 3, we can
prove $$\left|\displaystyle\frac{\partial(f\circ\phi)} {\partial
z_k}(z^{j_s})\right|(1-|z^{j_s}_k|^2)^q\to 0,$$ as $s\to\infty.$

Combining Case 1, Case 2, Case 3 and Case 4, we know there exists
a subsequence $\{z^{j_s}\}$ in $\{z^j\}$ such that
$$\left|\displaystyle\frac{\partial(f\circ\phi)} {\partial
z_k}(z^{j_s})\right|(1-|z^{j_s}_k|^2)^q\to 0,$$ as $s\to\infty$
for $k=1,2.$ We claim that
$$\left|\displaystyle\frac{\partial(f\circ\phi)} {\partial
z_k}(z^j)\right|(1-|z^j_k|^2)^q\to 0,$$ as $j\to\infty.$ In fact,
if it fails, then there exists a subsequence $\{z^{j_s}\}$ such
that
\begin{equation}\left|\displaystyle\frac{\partial(f\circ\phi)}
{\partial z_k}(z^{j_s})\right|(1-|z^{j_s}_k|^2)^q\to
\varepsilon>0\label{17}\end{equation} for $k=1$ or $2$. But from
the above discussion, we can find a subsequence in $\{z^{j_s}\}$
we still write $\{z^{j_s}\}$ with
$$\left|\displaystyle\frac{\partial(f\circ\phi)} {\partial
z_k}(z^{j_s})\right|(1-|z^{j_s}_k|^2)^q\to 0,$$  it contradicts
with (\ref{17}).

So for any sequence $\{z^j\}\subset U^n$ with $z^j\to\partial^*
U^n$ as $j\to\infty,$ we have
$$\left|\displaystyle\frac{\partial(f\circ\phi)} {\partial
z_k}(z^{j})\right|(1-|z^{j}_k|^2)^q\to 0$$ for $k=1,2.$ By
(\ref{9}) and Lemma 3, it is clear that $f\circ\phi\in{\cal
B}^q(U^n),$ so $f\circ\phi\in{\cal B}^q_{0*}(U^n).$

For any $f\in {\cal B}^p_0(U^n)).$ Since ${\cal
B}^p_0(U^n))\subset {\cal B}^p_{0*}(U^n)),$ then
$f\circ\phi\in{\cal B}^q_{0*}(U^n).$

By closed graph theorem we known that $$C_{\phi}:{\cal
B}^p_{0*}(U^{n}) ({\cal B}^p_0(U^n)) \rightarrow{\cal
B}^q_{0*}(U^{n})$$ is bounded.  This ends the proof of Lemma 4.

{\bf Remark 1}\hspace*{4mm}For the case $C_{\phi}:{\cal
B}^p(U^n)\to{\cal B}^q_{0*}(U^n)$, the necessity is also true, but
we can't guaranty that the sufficiency is true because we can't
sure that $C_{\phi}f\in{\cal B}^q_{0*}(U^n)$ for all $f\in{\cal
B}^p(U^n$.

\begin{Lemma} Let $\phi=(\phi_1, \phi_2, \cdots, \phi_n)$ be a holomorphic
self-map of $U^n.$ Then $$C_{\phi}:{\cal B}^p_0(U^n)
\rightarrow{\cal B}^q_0(U^{n})$$ is bounded if and only if if and
only if $\phi^{\gamma}\in {\cal B}_0^q(U^n)$ for every multi-index
$\gamma$, and (\ref{9}) holds.\end{Lemma}

{\bf Proof} \hspace{2mm}Sufficiency. From (\ref{9}) and by Theorem
1 we know that $C_{\phi}:{\cal B}^p(U^n)\to{\cal B}^q(U^n)$ is
bounded, in particular
$$\|C_\phi f\|_q\leq \|C_\phi\|_{{\cal B}^p(U^n)\to
{\cal B}^q(U^n)}\|f\|_p,\quad \mbox{for all}\; f\in{\cal
B}_0^p(U^n).$$ The boundedness of $C_{\phi}: {\cal B}_0^p(U^n)\to
{\cal B}_0^q(U^n)$ directly follows, if we prove $C_\phi f\in{\cal
B}_0^q(U^n)$ whenever $f\in {\cal B}_0^p(U^n).$ So, let $f\in
{\cal B}_0^p(U^n).$ By the definition of ${\cal B}_0^p(U^n)$ it
follows that for every $\varepsilon>0$ there is a polynomial
$p_\varepsilon$ such that $\|f-p_\varepsilon\|_p<\varepsilon.$
Hence \begin{equation}\|C_\phi f-C_\phi p_\varepsilon\|_q\leq
\|C_\phi\|_{{\cal B}^p(U^n)\to {\cal
B}^q(U^n)}\|f-p_\varepsilon\|_p<\varepsilon \|C_\phi\|_{{\cal
B}^p(U^n)\to {\cal B}^q(U^n)}.\label{a}\end{equation} Since
$\phi^{\gamma}\in{\cal B}_0^q(U^n)$ for every multi-index
$\gamma,$ we obtain $C_\phi p_\varepsilon\in{\cal B}_0^q(U^n).$
From this and (\ref{a}) the result follows.

If $C_{\phi}:{\cal B}_0^p(U^n)\to{\cal B}_0^q(U^n)$ is bounded,
then (\ref{9}) can be proved as in  Lemma 3, since the test
functions appearing there belong to ${\cal B}_0^p(U^n).$ Since the
polynomials $z^\gamma\in {\cal B}_0^p(U^n)$ for every multi-index
$\gamma,$ we get $C_\phi z^\gamma\in {\cal B}_0^q(U^n),$ as
desired.

{\bf Remark 2}\hspace*{4mm}For the case $C_{\phi}:{\cal
B}^p(U^n)\;\; ({\cal B}^p_{0*}(U^n))\to{\cal B}^q_{0}(U^n)$,
similar to Remark 1, the necessity is also true, but we can't
guaranty that the sufficiency is true.

\begin{Lemma} If $\{f_k\}$ is a bounded sequence in
${\cal B}^p(U^n)$, then there exists a subsequence $\{f_{k_l}\}$
of $\{f_k\}$ which converges uniformly on compact subsets of $U^n$
to a holomorphic function $f\in{\cal B}^p(U^n)$.
\end{Lemma}

{\bf Proof}\hspace{2mm} Let $\{f_k\}$ be a bounded sequence in
${\cal B}^p(U^n)$ with $\|f_k\|_p\leq C.$ By Lemma 1, $\{f_j\}$ is
uniformly bounded on compact subsets of $U^n$ and hence normal by
Montel's theorem. Hence we may extract subsequence $\{f_{j_k}\}$
which converges uniformly on compact subsects of $U^n$ to a
holomorphic function $f$. It follows that
$\displaystyle\frac{\partial f_{j_k}}{\partial
z_l}\to\displaystyle\frac{\partial f}{\partial z_l}$ for each
$l\in\{1,2,\cdots,n\}$, so
$$\sum\limits^n_{l=1}\left|\displaystyle\frac{\partial f}{\partial z_l}\right|(1-|z_l|^2)^p=
\lim\limits_{k\to\infty}\sum\limits^n_{l=1}\left|\displaystyle\frac{\partial
f_{j_k}}{\partial z_l}\right|(1-|z_l|^2)^p=\leq
\sup\limits_{k}\|f_{j_k}\|_p\leq C,$$ which implies $f\in{\cal
B}^p(U^n)$. The Lemma is proved.

\begin{Lemma} Let $\Omega$ be a domain in $\mbox{\Bbb C}^n,$
$f\in H(\Omega).$ If a compact set $K$ and its neighborhood $G$
satisfy $K\subset G\subset\subset \Omega$ and $\rho=dist(K,
\partial G)>0,$ then
$$\sup\limits_{z\in K}\left|\displaystyle\frac{\partial f}{\partial z_j}(z)
\right|\leq\displaystyle\frac{\sqrt{n}}{\rho}\sup\limits_{z\in G}|f(z)|.$$
\end{Lemma}

{\bf Proof}\hspace{2mm}Since $\rho=dist(K, \partial G)>0,$ for any $a\in K,$
the polydisc $$P_a=\left\{(z_1, \cdots, z_n)\in\mbox{\Bbb C}^n: |z_j-a_j|
<\displaystyle\frac{\rho}{\sqrt{n}}, j=1,\cdots,n\right\}$$
is contained in $G.$ By Cauchy's inequality,
$$\left|\displaystyle\frac{\partial f}{\partial z_j}(a)
\right|\leq\displaystyle\frac{\sqrt{n}}{\rho}
\sup\limits_{z\in\partial^* P_a}|f(z)|\leq
\displaystyle\frac{\sqrt{n}}{\rho}\sup\limits_{z\in G}|f(z)|.$$
Taking the supremum for $a$ over $K$ gives the desired inequality.

\section{The Proof of Theorem 1}
Now we turn to the proof of Theorem 1.

The lower estimate. It is clear that $\{m^{p-1}z^m_1\}\subset{\cal
B}^p_0(U^n) \subset{\cal B}_{0*}(U^n) \subset{\cal B}(U^n)$ for
$m=1,2,\cdots,$ and this sequence converges to zero uniformly on
compact subsets of the unit polydisc $U^n.$
\begin{equation}\|m^{p-1}z^m_1\|_p =\sup\limits_{z\in U^n}
(1-|z_1|^2)^p|m^pz^{m-1}_1|.\label{18}\end{equation} Let
$p(x)=m^p(1-x^2)^px^{m-1},$ then
$$p'(x)=-m^px^{m-2}(1-x^2)^{p-1}\left[(2p+m-1)x^2-(m-1)\right],$$ so
$p'(x)\leq 0$ for $x\in
\left[\sqrt{\displaystyle\frac{m-1}{2p+m-1}},1\right],$ and
$p'(x)\geq 0$ for $x\in
\left[0,\sqrt{\displaystyle\frac{m-1}{2p+m-1}}\right].$

That is, $p(x)$ is a decreasing function for $x\in
\left[\sqrt{\displaystyle\frac{m-1}{2p+m-1}},1\right]$ and $p(x)$
is a increasing function for $x\in
\left[0,\sqrt{\displaystyle\frac{m-1}{2p+m-1}}\right].$ Hence
$$\max\limits_{x\in
[0,1]}p(x)=p\left(\sqrt{\displaystyle\frac{m-1}{2p+m-1}}\right).$$
It follows from (\ref{18}) that
$$\|m^{p-1}z^m_1\|_p
=p\left(\sqrt{\displaystyle\frac{m-1}{2p+m-1}}\right)=\left(\displaystyle\frac{2p}{2p+m-1}\right)^pm^p
\left(\displaystyle\frac{m-1}{2p+m-1}\right)^{\frac{m-1}{2}}
\to\left(\displaystyle\frac{2p}{e}\right)^p,$$ as $m\to\infty.$

Therefore, the sequence $\{m^{p-1}z^m_1\}_{m\geq 2}$ is bounded
away from zero. Now we consider the normalized sequence
$\{f_m=\displaystyle\frac{m^{p-1}z^m_1}{\|m^{p-1}z^m_1\|_p}\}$
which also tends to zero uniformly on compact subsets of $U^n.$
For each $m\geq 2,$ we define $$A_m=\{z=(z_1, \ldots, z_n)\in U^n:
r_m\leq |z_1|\leq r_{m+1}\},$$ where
$r_m=\sqrt{\displaystyle\frac{m-1}{2p+m-1}}.$ So
\begin{eqnarray*}&&\min\limits_{A_m}
\sum\limits^n_{l=1} \left\{\left|\displaystyle\frac{\partial f_m}
{\partial z_l}(z)\right|(1-|z_l|^2)^p\right\}
=\min\limits_{A_m}\left|\displaystyle\frac{\partial f_m}
{\partial z_1}(1-|z_1|^2)^p\right|\\
&&=\displaystyle\frac{
(1-r^2_{m+1})^p|m^pr^{m-1}_{m+1}|}{\|m^{p-1}z^m_1\|_p}
=\left(\displaystyle\frac{2p+m-1}{2p+m}\right)
\left(\displaystyle\frac{m(2p+m-1)}{(m-1)(2p+m)}\right)^{\frac{m-1}{2}}=c_m.
\end{eqnarray*}

It is easy to show that $c_m$ tends to 1 as $m\to\infty$. For the
moment fix any compact operator $K:{\cal B}^p(U^n){\cal B}^p(U^n)
({\cal B}^p_0(U^n)$ or ${\cal B}^p_{0*}(U^n))\longrightarrow {\cal
B}^q(U^n)$ $({\cal B}^q_0(U^n)$ or ${\cal B}^q_{0*}(U^n)).$
 The
uniform convergence on compact subsets of the sequence $\{f_m\}$
to zero and the compactness of $K$ imply that $\|Kf_m\|_q\to 0.$
It is easy to show that if a bounded sequence that is contained in
${\cal B}^p_{0*}(U^n)$ converges uniformly on compact subsets of
$U^n,$ then it also converges weakly to zero in ${\cal
B}^p_{0*}(U^n)$ as well as in ${\cal B}^p(U^n).$ Since
$\|f_m\|_p=1$, we have
\begin{eqnarray*}&&\|C_{\phi}-K\|\geq
\limsup\limits_{m}\|(C_{\phi}-K)f_m\|_q\nonumber\\
&&\geq\limsup\limits_{m}\left(\|C_{\phi}f_m\|_q -\|Kf_m\|_q\right)
=\limsup\limits_{m}\|C_{\phi}f_m\|_q\nonumber\\
&&\geq \limsup\limits_{m}\sup\limits_{z\in U^n}\sum\limits^n_{k=1}
\left\{\left|\displaystyle\frac{\partial(f_m\circ\phi)}{\partial
z_k}\right| (1-|z_k|^2)^q
\right\}\nonumber\\
&&=\limsup\limits_{m}\sup\limits_{z\in U^n}\sum\limits^n_{k=1}
\left|\displaystyle\frac{\partial f_m}{\partial
w_1}(\phi(z))\right| \left|\displaystyle\frac
{\partial\phi_1}{\partial z_k}(z)\right|(1-|z_k|^2)^q\nonumber\\
&&=\limsup\limits_{m}\sup\limits_{z\in U^n}\sum\limits^n_{k=1}
\left|\displaystyle\frac{\partial\phi_1}{\partial z_k}(z)\right|
\displaystyle\frac{(1-|z_k|^2)^q}{(1-|\phi_1(z)|^2)^p}\left|\displaystyle\frac
{\partial f_m}{\partial w_1}(\phi(z))\right|(1-|\phi_1(z)|^2)^p\nonumber\\
&&\geq\limsup\limits_{m}\sup\limits_{\phi(z)\in A_m}
\sum\limits^n_{k=1}\left|
\displaystyle\frac{\partial\phi_1}{\partial
z_k}(z)\right|\displaystyle\frac
{(1-|z_k|^2)^q}{(1-|\phi_1(z)|^2)^p}\left|\displaystyle\frac
{\partial f_m}{\partial w_1}(\phi(z))\right|(1-|\phi_1(z)|^2)^p\nonumber\\
&&\geq\limsup\limits_{m}\sup\limits_{\phi(z)\in
A_m}\sum\limits^n_{k=1}\left|
\displaystyle\frac{\partial\phi_1}{\partial z_k}(z)\right|
\displaystyle\frac
{(1-|z_k|^2)^q}{(1-|\phi_1(z)|^2)^p}\nonumber\\
& &\times\liminf\limits_{m}\min\limits_{\phi(z)\in A_m}
\left|\displaystyle\frac{\partial f_m}
{\partial w_1}(\phi(z))\right|(1-|\phi_1(z)|^2)^p\nonumber\\
&&\geq\limsup\limits_{m} \sup\limits_{\phi(z)\in
A_m}\sum\limits^n_{k=1}\left|
\displaystyle\frac{\partial\phi_1}{\partial
z_k}(z)\right|\displaystyle\frac
{(1-|z_k|^2)^q}{(1-|\phi_1(z)|^2)^p}\liminf\limits_m c_m\nonumber\\
&&\geq\limsup\limits_{m} \sup\limits_{\phi(z)\in
A_m}\sum\limits^n_{k=1}\left|
\displaystyle\frac{\partial\phi_1}{\partial
z_k}(z)\right|\displaystyle\frac
{(1-|z_k|^2)^q}{(1-|\phi_1(z)|^2)^p}.
\end{eqnarray*}
So \begin{eqnarray}\|C_{\phi}\|_e
&=&\inf\{\|C_{\phi}-K\|: K \mbox{ is compact}\}
\nonumber\\
&\geq&\limsup\limits_{m} \sup\limits_{\phi(z)\in
A_m}\sum\limits^n_{k=1}\left|
\displaystyle\frac{\partial\phi_1}{\partial
z_k}(z)\right|\displaystyle\frac
{(1-|z_k|^2)^q}{(1-|\phi_1(z)|^2)^p}.\label{19}
\end{eqnarray}
For each $l=1,2,\cdots, n, $ define
\begin{equation}a_l=\lim\limits_{\delta\to 0}
\sup\limits_{dist(\phi(z), \partial U^n)<\delta}
\sum\limits^n_{k=1}\left|
\displaystyle\frac{\partial\phi_l}{\partial
z_k}(z)\right|\displaystyle\frac
{(1-|z_k|^2)^q}{(1-|\phi_l(z)|^2)^p}.\label{20}\end{equation} For
any $\varepsilon>0,$ (\ref{20}) shows that there exists a
$\delta_0$ with $0<\delta_0<1,$ such that
\begin{equation}\sum\limits^n_{k=1}\left|
\displaystyle\frac{\partial\phi_l}{\partial
z_k}(z)\right|\displaystyle\frac
{(1-|z_k|^2)^q}{(1-|\phi_l(z)|^2)^p}>a_l-\varepsilon,\label{21}\end{equation}
whenever $dist(\phi(z),\partial U^n)<\delta_0$ and
$l=1,2,\cdots,n.$ Since $r_m\to 1$ as $m\to\infty,$  so as $m$ is
large enough, $r_m>1-\delta_0.$ If $\phi(z)\in A_m,$ $r_m\leq
|\phi_1(z)|\leq r_{m+1},$ so
$1-r_{m+1}<1-|\phi_1(z)|<1-r_m<\delta_0,$ $dist(\phi_1(z),\partial
U)<\delta_0.$ There exists $w_1$ with $|w_1|=1$ such that
$dist(\phi_1(z),w_1)=dist(\phi_1(z),\partial U)<\delta_0.$  Let
$w=(w_1, \phi_2(z),\ldots, \phi_n(z)),$ $w\in\partial U^n$, then
$$dist(\phi(z),\partial
U)\leq dist(\phi(z), w)=dist(\phi_1(z),w_1)<\delta_0.$$ By
(\ref{21}), (\ref{19}) implies that
$$\|C_{\phi}\|_e\geq a_1-\varepsilon.$$
Similarly, if we choose
$g_m(z)=\displaystyle\frac{m^{p-1}z^{m}_l}{\|m^{p-1}z^m_l\|}$, we
have
$$\|C_{\phi}\|_e\geq a_l-\varepsilon,$$
for every $l=2\cdots,n.$ So
\begin{eqnarray*}\|C_{\phi}\|_e&\geq&\displaystyle\frac{1}{n}
\sum\limits^n_{l=1}
(a_l-\varepsilon)\\
&=&\displaystyle\frac{1}{n}\sum\limits^n_{l=1}
(\lim\limits_{\delta\to 0} \sup\limits_{dist(\phi(z), \partial
U^n)<\delta} \sum\limits^n_{k=1}\left|
\displaystyle\frac{\partial\phi_l}{\partial
z_k}(z)\right|\displaystyle\frac
{(1-|z_k|^2)^q}{(1-|\phi_l(z)|^2)^p}-\varepsilon)\\
&\geq&\displaystyle\frac{1}{n} \lim\limits_{\delta\to 0}
\sup\limits_{dist(\phi(z), \partial U^n)<\delta}
\sum\limits^n_{k,l=1}
\left|\displaystyle\frac{\partial\phi_l}{\partial z_k}(z)\right|
\displaystyle\frac{(1-|z_k|^2)^q}{(1-|\phi_l(z)|^2)^p}-
\varepsilon.\end{eqnarray*} Let $\varepsilon\to 0,$ the low
estimate follows.

The upper estimate. To obtain the upper estimate we first prove
the following proposition.

\begin{Proposition}  Let $\phi=(\phi_1, \ldots, \phi_n)$ a holomorphic self-map of $U^{n}.$ The operators
$K_m$ ($m\geq 2$) as follows:
$$K_mf(z)=f(\displaystyle\frac{m-1}{m}z),$$ for
$f\in H(U^n).$ Then the operators $K_m$ have the following
properties:

(i)\hspace*{2mm} For any $f\in H(U^n),$ $K_mf\in {\cal
B}^p_0(U^n)\subset {\cal B}^p_{0*}(U^n)\subset {\cal B}^p(U^n).$

(ii)\hspace*{2mm} If  $C_{\phi}:$ ${\cal B}^p(U^n)$ ( ${\cal
B}^p_{0}(U^n)$ or ${\cal B}^p_{0*}(U^n)$)$\rightarrow$ $ {\cal
B}^q(U^n)$ ( ${\cal B}^q_{0}(U^n)$ or ${\cal B}^q_{0*}(U^n)$) is
bounded, then $C_{\phi}K_mf\in{\cal B}^q(U^n)$ (${\cal
B}^q_{0}(U^n)$ or ${\cal B}^q_{0*}(U^n)$) for all $f\in H(U^n).$

(iii) \hspace*{2mm}For fixed $m$, the operator $K_m$ is compact on
${\cal B}^p(U^n)$ (${\cal B}^p_{0}(U^n)$ or ${\cal
B}^p_{0*}(U^n)$).

(iv)\hspace*{2mm} If $C_{\phi}:$ ${\cal B}^p(U^n)$ ( ${\cal
B}^p_{0}(U^n)$ or ${\cal B}^p_{0*}(U^n)$)$\rightarrow$ $ {\cal
B}^q(U^n)$ ( ${\cal B}^q_{0}(U^n)$ or ${\cal B}^q_{0*}(U^n)$) is
bounded, then $C_{\phi}K_mf\in{\cal B}^q(U^n)$ (${\cal
B}^q_{0}(U^n)$ or ${\cal B}^q_{0*}(U^n))$ is compact.

(v) \hspace*{2mm} $\|I-K_m\|\leq 2.$

(vi) \hspace*{2mm}$(I-K_m)f$ converges uniformly to zero on
compact subset of $U^n$.\end{Proposition}

{\bf Proof}\hspace*{2mm} (i)\hspace*{2mm} Let $f\in H(U^n),$
$r_m=\displaystyle\frac{m-1}{m}, (0<r_m<1)$ and
$f_m(z)=K_mf(z)=f(r_mz).$  First note that
\begin{eqnarray}\|f_m\|_p&=&|f(0)|+\sup\limits_{z\in U^n}\sum\limits^n_{k=1}
r_m\left|\frac{\partial f} {\partial z_k}(r_mz)\right|\left(1- |z_k|^2\right)^p\nonumber\\
&\leq&|f(0)|+\sup\limits_{z\in U^n}\sum\limits^n_{k=1}
\left|\frac{\partial f} {\partial z_k}(r_mz)\right|\left(1-
|r_mz_k|^2\right)^p\leq\|f\|_p.\label{b}\end{eqnarray}

On the other hand, $f_m\in H(\frac{1}{r_m}U^n).$
$0<\displaystyle\frac{2}{1+r_m}<\displaystyle\frac{1}{r_m},$
$\displaystyle\frac{2}{1+r_m}\overline{U^n}\subset
\displaystyle\frac{1}{r_m}U^n.$ which implies that for fixed $m,$
and $\varepsilon=\displaystyle\frac{1}{j}, j=1,2,\cdots,$ there is
a polynomial $P^{(j)}_m$ such that
$$\sup_{z\in \frac{2}{1+r_m}\overline{U^n}}|f_m(z)-P^{(j)}_m(z)|<(1-r_m)^2\displaystyle\frac{1}{j}.$$

Let $K=\overline{U^n},$ $G=\displaystyle\frac{2}{1+r_m}U^n,$
$\Omega=\displaystyle\frac{1}{r_m}U^n,$ then $K\subset
G\subset\subset\Omega$, and $\rho=dist(K,\partial
G)=\displaystyle\frac{1-r_m}{1+r_m}>0$, so $\forall w\in U^n$,
$k\in\{1,\cdots,n\}$, it follows from Lemma 7 that
\begin{eqnarray*}
&&\Big|\frac{\partial (f_m-P_m^{(j)})}{\partial w_k}(w)\Big|\leq
\sup_{w\in K}
\Big|\frac{\partial (f_m-P_m^{(j)})}{\partial w_k}(w)\Big|\\[6pt]
&\leq& \frac{\sqrt{n}(1+r_m)}{1-r_m}\sup_{w\in G}|f_m(w)-P_m^{(j)}(w)|\\[6pt]
&\leq& \frac{\sqrt{n}(1+r_m)}{1-r_m}(1-r_m^2)\frac{1}{j}\leq
4\sqrt{n}\frac{1}{j}.
\end{eqnarray*}
Therefore
$$
\sum_{k=1}^n\Big|\frac{\partial (f_m-P_m^{(j)})}{\partial
w_k}(w)\Big|(1-|w_k|^p)^p \leq 4n\sqrt{n}\frac{1}{j}\to 0
$$
as $j\to \infty.$ that is,
$$
||f_m-P_m^{(j)}||_{{\cal B}^p}=|f_m(0)-P_m^{(j)}(0)|+\sup_{w\in
U^n}\sum_{k=1}^n\Big|\frac{\partial (f_m-P_m^{(j)})}{\partial
w_k}(w)\Big|(1-|w_k|^p)^p\to 0.
$$
$P_m^{(j)}(w)\in {\cal B}^p_0(U^n)$ implies that $f_m\in{\cal
B}^p_0(U^n)$.

(ii)\hspace*{2mm} By (i), as desired.

(iii) \hspace{2mm} For any sequence $\{f_{j}\}\subset{\cal
B}^p(U^n)$ (${\cal B}^p_{0}(U^n)$ or ${\cal B}^p_{0*}(U^n)$) with
$\|f_j\|_p\leq M,$ by (i), $\{K_mf_{j}\}\in {\cal B}^p_{0}(U^n).$
By Lemma 6, there is a subsequence $\{f_{j_s}\}$ of $\{f_j\}$
which converges uniformly on compact subsets of $U^n$ to a
holomorphic function $f\in{\cal B}^p(U^n)$ and $\|f\|_p\leq M.$
$\left\{\displaystyle\frac{\partial f_{j_s}}{\partial z_i}
\right\}, i=1,2,\cdots,n,$ also converges uniformly on compact
subsets of $U^n$ to the holomorphic function
$\displaystyle\frac{\partial f} {\partial z_i}.$ So as $s$ is
large enough,  for any $w\in E=\{\frac{m-1}{m}z: z\in
\overline{U^n}\}\subset U^n$
\begin{equation}\left|\displaystyle\frac
{\partial (f_{j_s}-f)}{\partial
w_l}(w)\right|<\varepsilon,\label{23}\end{equation} for every
$l=1,2,\cdots,n.$ So
\begin{eqnarray}&&\left\|K_mf_{j_s}-K_mf\right\|_p
=\left\|f_{j_s}(\displaystyle\frac{m-1}{m}z)-
f(\displaystyle\frac{m-1}{m}z)\right\|_p\nonumber\\
&&=\sup\limits_{z\in U^n}\sum\limits^n_{k=1}
\left\{\left|\displaystyle\frac{\partial \left[(f_{j_s}-f)
(\displaystyle\frac{m-1}{m}z)\right]} {\partial z_k}\right|
(1-|z_k|^2)^p
\right\}+|f_{j_s}(0)-f(0)|\nonumber\\
&&\leq\sup\limits_{z\in
U^n}\sum\limits^n_{k=1}\sum\limits^n_{l=1}\left|
\displaystyle\frac{\partial (f_{j_s}-f)}{\partial w_l}
(\displaystyle\frac{m-1}{m}z)\right|\displaystyle\frac{m-1}{m}+|f_{j_s}(0)-f(0)|
\nonumber\\
&&\leq n\sup\limits_{w\in
E_1}\displaystyle\frac{m-1}{m}\sum\limits^n_{l=1}
\left|\displaystyle\frac {\partial (f_{j_s}-f)}{\partial
w_l}(w)\right|+|f_{j_s}(0)-f(0)|\to 0,\label{24}
\end{eqnarray}
as $s\to\infty.$ This shows that $\{K_mf_{j_s}\}$ converges to
$g=K_mf\in{\cal B}^p_{0}(U^n)\subset {\cal B}^p_{0*}(U^n)\subset
{\cal B}^p(U^n).$ So $K_m$ is compact on ${\cal B}^p(U^n)$(${\cal
B}^p_0(U^n)$ or ${\cal B}^p_{0*}(U^n)).$

(iv)\hspace*{2mm} By (i) and (iii), the result is obvious.

(v)\hspace*{2mm}In fact, for any $f\in{\cal B}^p(U^n)({\cal
B}^p_{0}(U^n)$ or ${\cal B}^p_{0*}(U^n))$, note that
$(I-K_m)f(0)=0$, so
\begin{eqnarray*}&&\|(I-K_m)f\|_p
=\sup\limits_{z\in U^n}\sum\limits^n_{k=1}
\left|\displaystyle\frac{\partial (I-K_m)f}{\partial z_k}(z)
\right|(1-|z_k|^2)\\
&&=n\sup\limits_{z\in U^n}\max\limits_{1\leq k\leq n}
\left|\displaystyle\frac{\partial f}{\partial z_k}(z)
-(1-\frac{1}{m}) \displaystyle\frac{\partial f}{\partial
z_k}((1-\frac{1}{m})z)
\right|(1-|z_k|^2)^p\\
&&\leq\sup\limits_{z\in U^n}\sum\limits^n_{k=1}
\left|\displaystyle\frac{\partial f}{\partial z_k}(z)\right|(1-|z_k|^2)^p\\
&&+(1-\frac{1}{m})\sup\limits_{z\in
U^n}\sum\limits^n_{k=1}\left|\displaystyle\frac{\partial
f}{\partial z_k}((1-\frac{1}{m})z)
\right|(1-|(1-\frac{1}{m})z_k|^2)^p\\
&&\leq \|f\|_p+\|f\|_p=2\|f\|_p,
\end{eqnarray*}
so $\|I-K_m\|\leq 2.$

(vi)\hspace*{2mm} For any compact subset $E\subset U^n$, $\exists
r,$ $0<r<1$ such that $E\subset rU^n\subset \subset U^n$. For
$\forall z\in E$,
\begin{eqnarray*}
|(I-K_m)f(z)|&=&|f(z)-f_m(z)|=|f(z)-f(r_mz)|\\
&=&\left|\int_{r_m}^1\frac{d}{dt}(f(tz))dt\right|=\left|\int_{r_m}^1\sum_{k=1}^n
\frac{\partial f}{\partial w_k}(tz)\cdot z_kdt\right|\\
&\leq& \sum_{k=1}^n\int_{r_m}^1\left|\frac{\partial f}{\partial
w_k}(tz)\right|dt.
\end{eqnarray*}
$t\in[r_m,1]$, $\forall z\in U^n,\hspace*{4mm}
|tz_k|=t|z_k|<|z_k|<r,$\hspace*{4mm} so
$\displaystyle\frac{\partial f}{\partial w_k}(w)$ is bounded in
$r\overline{U^n}$, i.e., $\forall z\in E,$\hspace*{4mm}
$\left|\displaystyle\frac{\partial f} {\partial
w_k}(tz)\right|\leq M$. So
$$
|(I-K_m)f(z)|\leq nM(1-r_m)\to 0$$ as $m\to \infty$, the results
follows.

Now return to the upper estimate. For the convenience, we denote
$\|f\|=\|f\|_p.$
\begin{eqnarray}&&\|C_{\phi}\|_e \leq\|C_{\phi}-C_{\phi}K_m\|
=\|C_{\phi}(I-K_m)\|
=\sup\limits_{\|f\|=1}\|C_{\phi}(I-K_m)f\|_q\nonumber\\
&&=\sup\limits_{\|f\|=1}\left(\sup\limits_{z\in
U^n}\sum\limits^n_{k=1} \left\{\left|\displaystyle\frac{\partial
(I-K_m)(f\circ\phi)} {\partial
z_k}\right|(1-|z_k|^2)^q\right\}+\left|(I-K_m)f(\phi(0))\right|\right)
\nonumber\\
&&\leq\sup\limits_{\|f\|=1}\sup\limits_{z\in
U^n}\sum\limits^n_{k=1}
\sum\limits^n_{l=1}\left|\displaystyle\frac{\partial (I-K_m)f}
{\partial w_l}(\phi(z))\right| \left|\displaystyle\frac
{\partial\phi_l}{\partial z_k}(z)\right|(1-|z_k|^2)^q\nonumber\\
&& +\sup\limits_{\|f\|=1}
\left|f(\phi(0))-f(\frac{m-1}{m}\phi(0))\right|\nonumber\\
&&\leq \sup\limits_{\|f\|=1}\sup\limits_{z\in U^n}
\sum\limits^n_{k, l=1}\left|
\displaystyle\frac{\partial\phi_l}{\partial z_k}(z)\right|
\displaystyle\frac {(1-|z_k|^2)^q}{(1-|\phi_l(z)|^2)^p}
\left|\displaystyle\frac{\partial (I-K_m)f}{\partial w_l}(\phi(z))
\right|(1-|\phi_l(z)|^2)^p\nonumber\\
&& +\sup\limits_{\|f\|=1}
\left|f(\phi(0))-f(\frac{m-1}{m}\phi(0))\right|.\label{26}\end{eqnarray}

Denote $G_1=\{z\in U^n: dist(\phi(z), \partial U^n)<\delta\},$
$G_2=\{z\in U^n: dist(\phi(z), \partial U^n)\geq\delta\},$
$G=\{w\in U^n: dist(w, \partial U^n)\geq\delta\}$, where $G$ is a
compact subset of $\mbox{\Bbb C}^n.$

Then by Lemma 3, Lemma 4 and Lemma 5, condition (9) holds, so
\begin{eqnarray}\|C_{\phi}\|_e
&\leq& \sup\limits_{\|f\|=1}\sup\limits_{z\in G_1}
\sum\limits^n_{k, l=1}\left|
\displaystyle\frac{\partial\phi_l}{\partial
z_k}(z)\right|\displaystyle\frac
{(1-|z_k|^2)^q}{(1-|\phi_l(z)|^2)^p} \left|\displaystyle\frac
{\partial (I-K_m)f}{\partial
w_l}(\phi(z))\right|(1-|\phi_l(z)|^2)^q \nonumber\\
&&+C\sup\limits_{\|f\|=1}\sup\limits_{z\in G_2}
\sum\limits^n_{l=1}(1-|\phi_l(z)|^2)^p
\left|\displaystyle\frac{\partial (I-K_m)f}{\partial
w_l}(\phi(z))\right|\nonumber\\
&& +\sup\limits_{\|f\|=1}
\left|f(\phi(0))-f(\frac{m-1}{m}\phi(0))\right|
\nonumber\\
&\leq& \|I-K_m\|\sup\limits_{z\in G_1} \sum\limits^n_{k,
l=1}\left| \displaystyle\frac{\partial\phi_l}{\partial
z_k}(z)\right|\displaystyle\frac
{(1-|z_k|^2)^q}{(1-|\phi_l(z)|^2)^p}\nonumber\\
&&+C\sup\limits_{\|f\|=1}\sup\limits_{z\in
G_2}\sum\limits^n_{l=1}(1-|\phi_l(z)|^2)^p
\left|\displaystyle\frac{\partial (I-K_m)f}{\partial w_l}(\phi(z))\right|\nonumber\\
&& +\sup\limits_{\|f\|=1}
\left|f(\phi(0))-f(\frac{m-1}{m}\phi(0))\right|\nonumber\\
&\leq& 2\sup\limits_{z\in G_1} \sum\limits^n_{k, l=1}\left|
\displaystyle\frac{\partial\phi_l}{\partial
z_k}(z)\right|\displaystyle\frac
{(1-|z_k|^2)^q}{(1-|\phi_l(z)|^2)^p}\nonumber\\
&&+C\sup\limits_{\|f\|=1}\sup\limits_{z\in G_2}
\sum\limits^n_{l=1}(1-|\phi_l(z)|^2)^p
\left|\displaystyle\frac{\partial (I-K_m)f}{\partial
w_l}(\phi(z))\right|\nonumber\\
&& +\sup\limits_{\|f\|=1}
\left|f(\phi(0))-f(\frac{m-1}{m}\phi(0))\right|.\label{27}
\end{eqnarray}
Denote the second term  and third term of the right hand side of
(\ref{27}) by $I_1$ and $I_2$. Then Theorem 1 is proved if we can
prove
$$\lim\limits_{m\to\infty}I_1=0\hspace*{4mm} \mbox{and}\hspace*{4mm} \lim\limits_{m\to\infty}I_2=0.$$
To do this, let $z\in G_2$ and $w=\phi(z),$ then $w\in G$
\begin{eqnarray}I_1&\leq&C\sup\limits_{\|f\|=1}
\sup\limits_{w\in G}\sum\limits^n_{l=1}(1-|w_l|^2)^p
\left|\displaystyle\frac{\partial f}{\partial w_l}(w)-
(1-\frac{1}{m})\displaystyle\frac{\partial f}{\partial w_l}
((1-\frac{1}{m})w)\right|\nonumber\\
&\leq& C\sup\limits_{\|f\|=1} \sup\limits_{w\in
G}\sum\limits^n_{l=1}(1-|w_l|^2)^p
\left|\displaystyle\frac{\partial f}{\partial w_l}(w)-
\displaystyle\frac{\partial f}{\partial w_l}
((1-\frac{1}{m})w)\right|\nonumber\\
& &+\displaystyle\frac{C}{m}\sup\limits_{\|f\|=1}
\sup\limits_{w\in G}\sum\limits^n_{l=1}(1-|w_l|^2)^p
\left|\displaystyle\frac{\partial f}{\partial w_l}
((1-\frac{1}{m})w)\right|\nonumber\\
&\leq& C\sup\limits_{\|f\|=1} \sup\limits_{w\in
G}\sum\limits^n_{l=1}(1-|w_l|^2)^p
\left|\displaystyle\frac{\partial f}{\partial w_l}(w)-
\displaystyle\frac{\partial f}{\partial w_l}
((1-\frac{1}{m})w)\right|+\displaystyle\frac{C}{m}.\label{28}
\end{eqnarray}
Let $w=(w_1,w_2,\cdots,w_{n-1},w_n),$ for $m$ large enough, we have
\begin{eqnarray}&&\left|\displaystyle\frac{\partial f}{\partial w_l}(w)-
\displaystyle\frac{\partial f}{\partial w_l}((1-\frac{1}{m})w)
\right|\nonumber\\
&&\leq\sum\limits^n_{j=1}\left|\displaystyle\frac{\partial f}
{\partial w_l}\left((1-\frac{1}{m})w_1,\cdots,
(1-\frac{1}{m})w_{j-1},
w_j,\cdots, w_n\right)\right.\nonumber\\
& &-\left.\displaystyle\frac{\partial f}{\partial w_l}
\left((1-\frac{1}{m})w_1,\cdots,(1-\frac{1}{m})
w_j,w_{j+1},\cdots,w_n\right)\right|
\nonumber\\
&&=\sum\limits^n_{j=1}\left|\int^{w_j}_{(1-\frac{1}{m})w_j}
\displaystyle\frac{\partial^2 f}{\partial w_l\partial w_j}
\left((1-\frac{1}{m})w_1,\cdots,
(1-\frac{1}{m})w_{j-1},\zeta,
w_{j+1},\cdots, w_n\right)d\zeta\right|\nonumber\\
&&\leq\frac{1}{m}\sum\limits^n_{j=1} \sup\limits_{w\in
G}\left|\displaystyle\frac {\partial^2 f}{\partial w_l\partial
w_j}(w)\right|.\label{29}\end{eqnarray} Denote $G_3=\left\{w\in
U^n:dist(w,\partial U^n)> \displaystyle\frac{\delta}{2}\right\},$
then $G\subset G_3\subset\subset U^n.$

Since $dist(G, \partial G_3)=\displaystyle\frac{\delta}{2},$ then
by Lemma 7, (\ref{29}) gives
\begin{equation}\left|\displaystyle\frac{\partial f}{\partial w_l}(w)-
\displaystyle\frac{\partial f}{\partial w_l}
((1-\frac{1}{m})w)\right| \leq\displaystyle\frac{2n\sqrt{n}}
{m\delta}\max\limits_{z\in G_3} \left|\displaystyle\frac{\partial
f}{\partial w_l}(w)\right|.\label{30}\end{equation} On the other
hand, on the unit ball of ${\cal B}^p(U^n)$, we have
$$\sup\limits_{z\in G_3}(1-|w_l|^2)^p\left|\displaystyle\frac{\partial f}
{\partial w_l}(w)\right|=\sup\limits_{dist(w,\partial
U^n)>\frac{\delta}{2}}
(1-|w_l|^2)^p\left|\displaystyle\frac{\partial f} {\partial
w_l}(w)\right|\leq \|f\|_p=1,$$ namely
\begin{equation}\sup\limits_{z\in G_3}\left|\displaystyle\frac{\partial f}
{\partial w_l}(w)\right|
\leq\displaystyle\frac{1}{1-\left(\frac{\delta}{2}\right)^2}=
\displaystyle\frac{4}{4-\delta^2}.\label{31}\end{equation}
Combining (\ref{28}), (\ref{30}) and (\ref{31}), imply
$$I_1\leq\displaystyle\frac{2n\sqrt{n}C}{m\delta}
\displaystyle\frac{4}{4-\delta^2}+\displaystyle\frac{C}{m}$$ and
$\lim\limits_{m\to\infty}I_1=0.$

Now we can prove $\lim\limits_{m\to\infty}I_2=0$. In fact,
\begin{eqnarray*}&&f(\phi(0))-f(\frac{m-1}{m}\phi(0))=
\int^{1}_{\frac{m-1}{m}}\displaystyle\frac{d f(t\phi(0))}{d
t}dt=\sum\limits^n_{l=1}\int^{1}_{\frac{m-1}{m}}
\phi_l(0)\displaystyle\frac{\partial
f}{\partial\zeta_l}(t\phi(0))dt.\end{eqnarray*} By Lemma 1, it
follows that for any compact subset $K\subset U^n$, $|f(z)|\leq
C_K \|f\|_p=C_K.$ Let $K=\{z\in U^n: |z_i|\leq|\phi_i(0)|\},$
So\begin{eqnarray*}&&|f(\phi(0))-f(\frac{m-1}{m}\phi(0))|\leq
\sum\limits^n_{l=1}|\phi_l(0)|\int^{1}_{\frac{m-1}{m}}C_K dt \leq
nC_K(1-\frac{m-1}{m})=\frac{nC_K}{m},
\end{eqnarray*} so $I_2\leq \frac{nC_K}{m}\to 0.$
Thus let first $m\to\infty,$ then $\delta\to 0$ in (\ref{27}), we
get the upper estimate of $\|C_{\phi}\|_e$:
$$\|C_{\phi}\|_e\leq 2\lim\limits_{\delta\to 0}
\sup\limits_{dist(\phi(z),\partial
U^n)<\delta}\sum\limits^n_{k,l=1}
\left|\displaystyle\frac{\partial \phi_{l}}{\partial
z_k}(z)\right|
\displaystyle\frac{(1-|z_k|^2)^q}{(1-|\phi_l(z)|^2)^p}.$$ Now the
proof of Theorem 1 is finished.

\end{document}